%% file: arxinv.tex
\documentclass[times,doublespace]{amsart}

\usepackage[mathscr]{eucal}
\usepackage{amssymb, amsmath,array, amscd}
\usepackage{enumerate}
\usepackage[colorlinks]{hyperref}
\usepackage[all]{xy}

\def\rev{\usepackage[active]{srcltx}}
\def\p#1{\ensuremath{\mathbb{P}^{#1}}}
\def\pd#1{\ensuremath{\check{\mathbb{P}}^{#1}}}
\def\eg{{\em e.g., }}
\def\ie{{\em i.e., }}
\def\C{\ensuremath{\mathbb C}}
\def\ve{\ensuremath{^\vee}}
\def\lar{\ensuremath{\,\longrightarrow\,}}
\def\be{\begin{equation}}\def\ee{\end{equation}}\def\bsm{\left(\begin{smallmatrix}}\def\esm{\end{smallmatrix}\right)}
\def\ba#1{\begin{array}{#1}}\def\ea{\end{array}}
\def\cl#1{\ensuremath{\mathcal{#1}}}

\def\vez {\ensuremath{\times}}
\def\ra{\ensuremath{\,\rightarrow\,}}
\def\us{\ensuremath{^\star}}

\def\ds#1{\displaystyle{#1}}
\def\M{{\cl  M}}
\def\N{{\cl  N}}
\def\X{{\cl  X}}

\def\J{{\cl  J}}
\def\na#1{\noalign{\vskip#1pt}}
\def\wt#1{\ensuremath{\widetilde{#1}}}

\newtheorem{coisa}{}[section]
\def\coi#1{\begin{coisa}{\bf
      #1.} \em}
\def\eco{\end{coisa}}
\def\cois#1{\begin{coisa}{\bf #1.}
\vspace{-10pt}
\eco
\addcontentsline{toc}{subsection}
{\numberline{}{#1}}
}

\def\runningheads#1#2{\markboth{\MakeTextUppercase{#1}}{\MakeTextUppercase{#2}}}

\def\rk{\ensuremath{\operatorname{rank}}}
\def\Hom{\ensuremath{\operatorname{Hom}}}
\def\cod{\ensuremath{\operatorname{cod}}}

\def\V{\ensuremath{\mathbb{V}}}
\def\ide#1{\ensuremath{\langle{#1}\rangle}}
\def\s#1{\ensuremath{\operatorname{Sym}_{#1}}}
\def\wed#1{\ensuremath{\stackrel{#1}\wedge}}

\def\E{\ensuremath{\mathcal{E}}}
\def\schub{{\sc schubert}\,\cite{schub}}

\newcommand{\Y}{\ensuremath{\mathbb Y}}

\newcommand{\graf}{{\ensuremath{\mathrm{Graph}\,}}}
\newcommand{\Z}{{\mathbb Z}}
\newcommand{\grass}{{\mathbb G}}
\newcommand{\CC}{{\mathcal C}}
\newcommand{\de}{\partial}
\newcommand{\OX}{\mathscr{O}}

\newcommand{\T}{{\mathcal T}}
\newcommand{\Ker}{{\ensuremath{\mathrm{ker}}}}

\newcommand{\im}{{\ensuremath{\mathrm{Im}}}}

\def\inv{\ensuremath{{}^{-1}}}
\def\rref#1{(\ref{#1},\,p.\,\pageref{#1})}
\def\z{\ensuremath{\check z}}
\newcommand{\B}{\ensuremath{\mathbb B}}
\def\sing{{\sc Singular}\,\cite{sing}}
\def\mn{\medskip\noindent}
\rev
\begin{document}

\input{invmain}
\section{scripts}
\label{scripts}
\input{scripts1}

\input{scripts2}

\input{biblioinv}

%% file: invmain.tex
\runningheads{V. Ferrer and I. Vainsencher}
{vector fields  with algebraic trajectories}
\title[Polynomial vector fields  with algebraic trajectories]{Polynomial vector fields  \\with algebraic trajectories}
\author{V. Ferrer
and I. Vainsencher}

\address{
IMPA --- Estrada Dona Castorina 110. 224620-320 Rio de
Janeiro Brazil. 
\newline\indent
ICEX--Depto. Matem\'{a}tica, UFMG -- Av. Antonio
Carlos 6627. 31270-901 Belo Horizonte Brazil. 
}

\thanks{The authors were partially supported by  CNPQ.}


\begin{abstract}
It is known after Jouanolou that a general
holomorphic  foliation of degree
$\geq2$ in projective space has no algebraic
leaf.  We give formulas for the
degrees of the subvarieties of the parameter
space of one-dimensional foliations that
correspond to foliations endowed with some
invariant subvariety of degree $1$ or $2$ and
dimension $\geq1$.

\end{abstract}


\maketitle

\section*{Introduction}
Holomorphic foliations are an offspring of the
geometric theory of polynomial differential 
equations. Following the trend of
many branches in Mathematics,
interest has migrated to global aspects.
Instead of focusing on just one curve, or
surface, or metric, or differential equation,
try and study their family in a suitable
parameter space. The geometry within the
parameter space of the family acquires
relevance.  For instance, the family of
hypersurfaces of a given degree correspond to
points in a suitable projective space;
geometric conditions imposed on hypersurfaces
, \eg to be singular, 
usually correspond to interesting subvarieties
in the parameter space, \eg the discriminant.
Hilbert schemes have their counterpart in the
theory of polynomial differential equations, to
wit, the spaces of foliations.

Let us recall that, while a general surface of
degree $d\geq4$ in 
\p3 contains no line --in fact, only complete
intersection curves are allowed, those that
{\em do contain} some line correspond to a
subvariety of codimension $d-3 $ and degree
$ \binom{d+1}{4}
(3d^4+6d^3+17d^2+22d+24)/4!  $ in a suitable \p
N.

Similarly, while a general holomorphic foliation, say in
\p2, of degree $d\geq2$ has no algebraic leaf,
those that {\em do have}, say an invariant
line, correspond to a subvariety of codimension
$d-1$ and degree $3
\binom{d+3}{4}$ in a suitable \p N.

Our goal  is to
give similar
closed formulas for the degrees of the
subvarieties of the parameter space of
one-dimensional foliations of degree $d$ on \p
n that correspond to foliations endowed with
some invariant subvariety of degree $1$ or
$2$ and dimension $\geq1$.
Imposing a linear invariant subvariety is 
easy, essentially due to the absence of
degenerations. The classical spaces of complete
quadrics help us handle the quadratic case. 
For higher degree, we don't know the answer.

\label{section Introduction}

\section{The space of foliations}
The main reference for this matterial is
Jouanolou\,\cite{Joua}. We call a {\em vector
field of degree} $d$ in $\p n$ a global section
of $T\p n\otimes \OX_{\p n}(d-1)$, for some
$d\geq 0$.

Denote by $S_d$ the space $ {\rm
Sym}_d(\C^{n+1})\ve $ of homogeneous
polynomials of degree $d$ in the variables
$z_0,\dots,z_n$.  We write
$\partial_i=\partial/ \partial z_i $, thought
of as a vector field basis for
$\C^{n+1}=S_{1}^{\vee}$.  Recalling Euler
sequence
\begin{equation*}
 0\to \OX_{\p n}(d-1)\lar 
\OX_{\p n}(d)\times \C^{n+1}\lar
T\p n(d-1)\to 0
\end{equation*}  and taking global sections we get the exact sequence
\begin{equation}\label{globalEuler}
0\to S_{d-1}\stackrel{\iota}\lar 
S_{d}\otimes S_{1}^{\vee}\lar{}
V_d^n:=H^0(\p n,T\p n(d-1))\to 0. 
\end{equation} 
Here $\iota(H)=HR$, with $R=\sum z_i\partial_i$
the radial vector field. Any degree $d$ vector
field $\X\in V_d^n$ can be written in
homogeneous coordinates as
\be\label{X}
\X=F_{0}\partial_0+
\dots+F_{n}\partial_n ,
\ee 
where the $F_{i}$'s denote homogeneous
polynomials of degree $d$, modulo multiples of
the radial vector field.
A vector field $\X$ induces a distribution of
directions in $T\p n$. Any nonzero multiple of
$\X$ yields the same distribution.  

The {\em
space 
of foliations of degree} $d$ in  $\p n$ is the
projective space
\be\label{H0}
\ba c
\p N
=\p{}(V_d^n)
\ea\ee
of dimension

\centerline{$
N=(n+1)\binom{d+n}n-\binom{d+n-1}n-1
$.}  

\mn
We shall often abuse notation and denote by the
same symbol \X\ both the foliation and a vector
field. The {\em singular scheme} 
of \X\ is the scheme of zeros of the section
$\X:\OX\ra\T\p n(d-1)$. With \X\ as in (\ref{X}),
the singular
scheme is given by the 2\vez2 minors of the
matrix 
\\\centerline{$
\left(\ba{cccc}
F_0&F_1&\cdots &F_n\\
z_0&z_1&\cdots&z_n
\ea\right).$}

\coi{Invariant subvarieties}
Let $\X$ be a foliation in $\p n$. An
irreducible subvariety $Z\subset \p n$ is said
to be {\em invariant by} $\X$ if
$$ 
\X(p)\in T_{p}Z
$$ 
for all $p\in Z\setminus (Sing(Z)\cup
Sing(\X))$. 
If $Z$ is reducible,  it is invariant by $\X$ if and only if each irreducible component of $Z$ is invariant by $\X$.
If $Z$ is defined  by a saturated ideal
$I_{Z}:=\langle G_{1},\dots,G_{r}\rangle$,
invariance means  
$$dG_{i}(\X)=\X(G_{i})\in I_{Z}$$ for all
$i=1,\dots,r$.  
The hypothesis of saturation is
necessary, see \cite[p.\,5]{Est}. 
It can be easily checked that the condition above does not depend on the representative
of $\X$ in $S_{d}\otimes S_{1}^{\vee}$. 

\eco
\section{Foliations with an invariant 
$k$-plane}

We show that the locus in \p N corresponding to
foliations with an invariant $k$-plane is 
the birational image of a natural projective
bundle over the grassmannian of $k$-planes in
\p n.

\begin{coisa}
\em
Consider  the tautological exact sequence of
vector bundles over the grassmannian  
$\grass:=\grass(k,n)$ of projective $k$-planes 
in \p n,
\begin{equation}\label{tau}
0\to \T\to \grass\times\C^{n+1}\to \mathcal{Q}\to 0
\end{equation}
where $\T$ is of rank $k+1$. The
projectivization
$$\p {}(\T)=\{(W,p)\in \grass\times \p n\mid
p\in W\}
$$
is the {\em universal $k$-plane}.  Write the
projection maps
\[
\xymatrix 
{
&~
\p {}(\T) \subset\grass\vez\p n
\ar[dl]_{p_{1}} \ar[dr]^{p_{2}}&  \\
\grass & &\p n.}
\]
We denote by $\cl N$ the normal bundle to $\p
{}(\T)$ in $\grass\times \p n$.  We have the
exact sequence over $\p {}(\T)$,
\begin{equation}\label{relnormal}
0\rightarrow T_{\p {}(\T)/\grass} \rightarrow p_{2}^*T\p n
\rightarrow \cl N\rightarrow 0.
\end{equation}
It is easy to see that 
$$
\cl N =p_{1}^*\mathcal{Q}\otimes \OX_{\T}(1).
$$
\end{coisa}

\begin{coisa}{\bf Proposition.}
Notation as in {\rm(\ref{globalEuler})}
 and {\rm(\ref{tau})},
there exists a vector subbundle 

\medskip
\centerline{$\E\subset{}
\grass\vez V_d^n$}
\mn
such that 

\medskip
\centerline{$\ba l
{\rm(i)}\ \p{}(\E)=\{(W,\X)
\in  \grass\times \p N
\mid W\,\text{ is invariant by }\,\X\}.
\\\text{Set }
\Y:=q_2(\p{}(\E)), \text{ where }
q_2:\p{}(\E)\to \p N  \text{ is the projection.
Then} 
\\{\rm(ii)}  \text{ the
codimension of  $\Y$ in \p N is
$(n-k)(\binom{k+d}{d}-(k+1))$ and} 
\\{\rm(iii)}  \text{ 
the  degree of  $\Y$ is given by 
the top-dimensional Chern class,}
\ea$} 

\centerline{$
c_{g}(\mathcal{Q}\otimes \s {d}(\T^{\vee}))
$}

\noindent
where $g:=\dim \grass$.

\end{coisa}

\begin{proof}

Consider the following diagram of maps of 
vector bundles over $\p {}(\T)$,
\\\centerline{$
\hskip-1cm
\xymatrix@R1.95pc
{
     &    &T_{\p {}(\T)/\grass}(d-1)
\vphantom{\ba c1\\.\ea}
\ar@{>->}[d]\\
 &  V_d^n=H^0(\p n,T\p n(d-1))
\ar@{->}[r]^{\hskip.81cm
   \,ev} \ar@{->}[dr]_\varphi&  p_{2}^*T\p
 n(d-1) 
\ar@{->>}[d]\\
&   &    
\ \ \ \, p_{1}^*\mathcal{Q}\otimes \OX_{\T}(d).
}$}

\mn
The map of evaluation yields a surjective map
of vector bundles  
$$\varphi:
V_d^n
\twoheadrightarrow
p_{1}^*\mathcal{Q}\otimes \OX_{\T}(d).
$$
The  kernel of \,$\varphi$\, is the vector bundle
\[
\Ker\varphi=
\{(W,p,\X)\mid\, \X(p)\in T_{p}W\}.\]
Observe that a $k$-plane $W$ is invariant by
$\X$ if and only if $(W,p,\X)\in\Ker\varphi$ for all
$p\in W=p_{1}^{-1}(W)$. Put in other words, we
ask $\varphi$ to vanish along the fibers of
$p_1:\p{}(\T)\ra\grass$. It follows from
 \cite[p.\,16]{AK} that 
there exists a map of vector bundles over $\grass$,

\medskip\noindent
$(\heartsuit)
~\hfill
$\xymatrix
@C3pc
{\varphi^{\flat}:
V_d^n
\,\ar@{->>}[r]&
p_{1*}(p_{1}^*\mathcal{Q}\otimes
\OX_{\T}(d))=\mathcal{Q}\otimes \s
{d}(\T^{\vee})}$
\hfill
$

\medskip\noindent
such that 

\centerline{$
\E:=\Ker(\varphi^{\flat})=\{(W,\X)\mid\, \X(p)\in T_{p}W\,\,\forall p\in
W\}.$}

\medskip\noindent
The projectivization $
\p{}(\E)\subset\grass\vez\p N$ is clearly 
as stated in (i). Let $q_1:\p{}(\E)\ra
\grass$ and $q_2:\p{}(\E)\ra\p N$ be the 
maps induced by projection. It can be shown 
that $q_2:\p{}(\E)\ra\p N$ is
generically injective: the general vector field
of degree $d\geq2$ with an invariant $k$-plane
has exactly one  
invariant $k$-plane. Write  $H$ for the
hyperplane class of \p N. Set
$u=\dim\Y=\dim\p{}(\E)$. We have
$q_2^*H=c_1\OX_{\E}(1)=:h$. We may compute
$$
\deg\Y=\int_{\p N}H^u\cap\Y=\int_{\p{}(\E)}
h^u=\int_\grass q_{1*}(h^u)=\int_\grass{}
s_{g}(\E).
$$
The assertions (ii) and (iii) follow noticing
that  the Segre class satisfies
\\\centerline{$
s_{g}(\E)=c_{g}(\mathcal{Q}\otimes \s
{d}(\T^{\vee}))
$} 
in view of the exact sequence
arising from $(\heartsuit)$,

  \centerline{$
\xymatrix@C2.753pc
{\E\ \ 
\ar@{>->}[r]&
V_d^n
\ \ar@{->>}[r]& 
\mathcal{Q}\otimes \s{d}(\T^{\vee}).
}
$}

\end{proof}

In the case of invariant hyperplane 
(\ie $k=n-1$) we have explicitly
\begin{coisa}{\bf Theorem.}
The degree of the subvariety \Y\ of the space
of foliations in \p n that admit an invariant
hyperplane is given by
\\\centerline{$
\deg\Y=\ds{\binom{\binom{d+n}{n}}{n}}.
$}
~\hfill\qed
\end{coisa}

\begin{proof}
We have the following exact sequence
over $\grass(n-1,n)=\pd n$,
$$0\to S_{d-1}\otimes \mathcal{Q}^{\vee}\to 
S_{d}=\s d(\C^{n+1})\ve\to
\s{d}(\T^{\vee})\to 0.
$$
Twisting by $\mathcal{Q}=\OX_{\pd n}(1)$ we obtain:
$$0\to S_{d-1}\to S_{d}\otimes \mathcal{Q}\to
\s{d}(\T^{\vee})\otimes \mathcal{Q}\to 0.
$$
From this we can compute 
$c_{n}(\s{d}(\T^{\vee})\otimes
\mathcal{Q})=c_{n}(S_{d}\otimes \mathcal{Q})$.
Setting $H=c_{1}(\mathcal{Q})$, the hyperplane
class in \pd n, the sought for
degree is just  the coefficient of $H^n$ in
$(1+H)^{\binom{d+n}{n}}$.
\end{proof}

To compute 
$c_{g}(\mathcal{Q}\otimes \s {d}(\T^{\vee}))$
for any fixed $k,n$,  see the script in
\S\,\ref{scripts}.





\section{Foliations with an invariant conic in \p2}\label{conics}

We construct a compactification of the space of foliations that
 leave invariant a smooth conic. This compactification is obtained as
 the birational image of a projective bundle over the variety of
 complete conics.

\coi{The incidence variety}
The parameter space for the family of conics
is $\p 5=\p {}(S_{2})$. We have the natural
trilinear map
\\\centerline{$\ba c
S_2\otimes S_{d}\otimes S_{1}^{\vee}\to
S_{d+1}\\
G\otimes{}F\otimes{}\partial_i
\mapsto
(\partial_iG)F.
\ea$}
It induces  the 
map of vector bundles over \p5,
\\\centerline{$
\varphi:\OX_{\p 5}(-1)\otimes S_{d}\otimes S_{1}^{\vee}\to
S_{d+1}$}
given by 
$
\varphi(G,\X):=\X(G)=\sum
F_i\partial_iG
$, 
where $G$ is the equation of the
conic and $\X=\sum
F_i\partial_i,\,F_i\in{}S_d$.

Notice that $\varphi(G,HR)=2GH$, for all $H\in
S_{d-1}$. Recalling (\ref{globalEuler}),

\medskip
\centerline{$V_d^n
=S_{d}\otimes S_{1}^{\vee}\big/S_{d-1}R
$}
\noindent
we see that $\varphi$ induces a map
$$
\psi:\OX_{\p 5}(-1)\otimes V_d^n
\lar\frac{S_{d+1}}{\OX_{\p 5}(-1)\otimes S_{d-1}}
\,\cdot
$$
This maps fits into the  commutative
diagram,
\be\label{fipsi}
\xymatrix@C2.2pc
{&\Ker\varphi \ar@{->}[r]^{\simeq} \ar@{->}[d]&
 \ \Ker\psi  \ar@{->}[d]\\
\OX_{\p 5}(-1)\otimes S_{d-1}\cdot{}R\ 
\ar@{>->}[r] 
\ar@{->}[d]^{\simeq}& \OX_{\p 5}(-1)
\otimes{}S_{d}\otimes{}S_{1}^{\vee} 
\,\ar@{->>}[r]
\ar@{->}[d]^{\varphi}&\OX_{\p 5}(-1)\otimes{}
V_d^n \ar@{->}[d]^{\psi}  \\
\OX_{\p 5}(-1)\otimes S_{d-1}\ \ar@{>->}[r]&  
S_{d+1} \ar@{->>}[r]&  \,\frac{S_{d+1}}
{\OX_{\p 5}(-1)\otimes S_{d-1}}\cdot
}
\ee
Twisting by $\OX_{\p 5}(1)$ we obtain
\be\label{N}\xymatrix
{&\M \ar@{->}[r]^{\simeq} \ar@{->}[d]&  
\wt{\M} \ar@{->}[d]\\
S_{d-1}\cdot{}R \ \,
\ar@{>->}[r] \ar@{->}[d]^{\simeq}&
S_{d}\otimes S_{1}^{\vee}
\ar@{->>}[r]\ar@{->}[d]_\theta
& V_d^n
\ar@{->}[d]^{\wt\theta}
 \\
S_{d-1}\ \,
\ar@{>->}[r]& \OX_{\p 5}(1)\otimes S_{d+1}
\ar@{->>}[r]&  \,\frac{\OX_{\p 5}(1)\otimes{}
S_{d+1}}{S_{d-1}},
}
\ee
where $$\M:=
\Ker\,\theta\simeq\wt{\M}:=
\Ker\,\wt\theta.
$$

\eco

The map $\theta:S_{d}\otimes S_{1}^{\vee}\ra
\OX_{\p5}(1)\otimes{}S_{d+1} $ appearing in
(\ref{N}) is surjective only over the open
subset consisting of smooth
conics. 
In fact, for  $G\in\p5$
the rank of the image of $\theta_G$ depends
on the singularities of the conic:
$$\rk\theta_G= \left\{\ba l
\text{total}, 
\binom{d+3}{2}, \text{ if }G
\text{ is smooth;} \\
\binom{d+3}{2}-1, \text{ if }G
\text{ is a line pair;} \\
{\rm rank}\s d=\binom{d+2}2,
\text{ if }G
\text{ is a double line.}
\ea\right. $$
The minimal 
rank, $r=\binom{d+2}2$, is achieved along 
$$\V:=\nu_{2}(\p {}(S_{1}))\subset \p {}(S_{2})
,$$ 
the Veronese variety of double lines. 
It turns out that over the open subset 
\\\centerline{$U\subset \p5\setminus\V$}
of smooth conics, the restriction $\wt{\M}_{U}$
is a vector subbundle of the trivial bundle 
$V_d^n$.
The projectivization $\p{}(\wt{\M}_{U})$
is the 
incidence variety 
$$\{(\CC,\X)\mid \,\CC\, \text{is invariant
 by}\,\X\,\}\subset{}U\times \p {N}.
$$
Let us denote by \Y\ the closure of its image
in \p N. We see that \Y\
consists of (limits
of) foliations that admit an invariant  smooth
conic. Our strategy to find its degree is
summarized in the following.

\coi{Theorem}\label{invconics}
{\em 
Let $\pi:\mathbb{B}\to \p 5$ be the blowup of
$\p 5$ along the veronese $\V{}$. 
There exists a subbundle 
$$\E\subset \B\vez{}V_d^n
$$ 
such that the restriction
$\E_{\pi\inv U}$  coincides
with $\pi^*\wt{\M}_U$ \,{\rm (cf.\,\ref{N})}.
In particular, the image of $\p{}(\E)$ in 
$\p N=\p{}(V_d^n)$
is equal to \Y.
}

\eco

\begin{proof}

Consider the pullback by $\pi$ of the maps
$\varphi$ and $\psi$ (cf.\,\ref{fipsi})
$$\ba c
\varphi_{\mathbb{B}}:
\pi^*(\OX_{\p 5}(-1)\otimes{}S_{d}
\otimes{} S_{1}^{\vee})\to \pi^*S_{d+1},
\\\na7
\psi_{\mathbb{B}}:\pi^*(\OX_{\p 5}(-1)\otimes 
V_d^n
\to
\pi^*(\frac{S_{d+1}}{S_{d-1}\otimes \OX_{\p
    5}(-1)})\cdot
\ea
$$
By Lemma~\ref{locallema} below it's enough to prove that
 the $k\times k$ minors of
$\varphi_{\B}$ are locally principal for all $k\geq 1$. 
Indeed, in
this case 
 $\cl J:=\im\varphi_{\mathbb{B}}$
is locally free. Therefore we  obtain a
factorization of 
$\varphi_{\mathbb{B}}$,
\\\centerline{$
\xymatrix@R1.4pc@C3pc
{
\pi^*(\OX_{\p 5}(-1)\otimes
S_{d}\otimes S_{1}^{\vee}) \ar@{->}[r]
^{\hskip1.31cm\varphi_{\mathbb{B}}}
\ar@{->>}[dr]_{\ \widetilde{\varphi}}& 
\pi^*S_{d+1}\ar@{<-^)}[d]\\
           & 
\, \cl J\vphantom{\ba c.\\.\ea}
,
}$}
and $\pi^*(\OX_{\p 5}(-1)\otimes
S_{d-1})=\varphi_{\B}(\pi^*(\OX_{\p 5}(-1)\otimes
S_{d-1}\cdot{}R))\subset \J$.
This factorization induces a factorization of
$\psi_{\B}$ 
\[
\xymatrix@C3pc
{\pi^*(\OX_{\p 5}(-1)\otimes 
V_d^n
\ar@{->}[r]
^{\hskip-.02cm\psi_{\B}}
\ar@{->>}[dr]_{\widetilde{\psi}}& 
\pi^*(\frac{S_{d+1}}{\OX_{\p 5}(-1)\otimes
S_{d-1}})\ar@{<-}[d]^{\iota}\\
           &  \,\wt{\J},
}\]
where 
\\\centerline{$
\wt{\J}:=\ds{\frac{\J}{
 \OX_{\p 5}(-1)\otimes
\pi^*(S_{d-1}
)}}
$} 
is a vector bundle. Define
\be\label{E2}
\E:=\pi^*\OX_{\p 5}(1)\otimes
\Ker\widetilde{\psi}.
\ee
It follows that $\E$ is a subbundle of
\,$\B\times V_d^n$\, that coincides
with $\pi^*\wt{\M}$ over the open set
$\pi^{-1}(U)$, where $U\subset \p 5$ is the
open set of smooth 
conics. Indeed, over $\pi^{-1}(U)=U$ the map
\\\centerline{$
\iota:\wt{\J}\to
\pi^*\left(\ds{\frac{S_{d+1}}{\OX_{\p5}(-1)
\otimes{} S_{d-1}}}\right)
$}  

\medskip
\noindent
is an isomorphism. Therefore 
$\Ker\widetilde{\psi}=\Ker\psi_{\mathbb{B}}$
over $U$. 

To prove that the $k\times k$ minors of
a local matrix representation of
$\varphi_{\B}$ are  principal for all
$k\geq 1$,  
consider
$$
\varphi_{0}:\OX_{\p 5}(-1)\otimes
S_{1}^{\vee}\to S_{1},
$$ the universal
symmetric map that 
gives the matrix of the conic.
We are blowing-up the Veronese, which is the
scheme of zeros of  the ideal of $2\times
2$-minors  of 
$\varphi_{0}$. Therefore, up on \B,
we have that the 2\vez2 
minors of $\varphi_{0\B}$
are locally principal, say generated by $t$.
Thus we can write the matrix of
$\varphi_{0}$ locally in the form $A=\bsm
1&0&0\\
0&t&0\\
0&0&ts
\esm$ (cf.\,\cite{Laksov}). 
Let $\CC=Z(z_{0}^2+tz_{1}^2+tsz_{2}^2)$
be the associated conic.

Choosing
an appropriate ordering of the basis of
$S_{d}\otimes S_{1}^{\vee}$ and of $S_{d+1}$
the matrix of $\varphi_{\CC}$ can be put in
 the following form:
$$A_{d}=\left(\begin{array}{ccccc}
2I_{n(d)}&0&0&B_{1}&B_{3}\\
0&2tI_{d+1}&0&0&B_{4}\\
0&0&2ts&0&0
\end{array}\right),
$$
where all the entries of $B_{1}$ are  multiples
of $t$, and the 
entries of $B_{3},B_{4}$ are  multiples of $ts$.
Here $I_{m}$ stands for the identity matrix of
size $m$, and $n(d)=\binom{d+2}{2}$.
It follows that the ideals $J_{i}$ of $i\times
i$-minors of $A_{d}$ are:
\\\centerline{$
\left\{\ba{lcl}
J_{i}&=&
\langle 1\rangle\, \text{ for }\,i=1,\dots,n(d)
;\\
J_{n(d)+j}&=&\langle t^j\rangle\, \text{ for }
\, j=1,\dots,d+1;\\
J_{n(d+1)}&=&\langle t^{d+2}s\rangle.
\ea\right.$}
In particular these minors are principal as  
claimed.
\end{proof}

\begin{coisa}{\bf Lemma.}\label{locallema}
Let $R$ be a local Noetherian domain, and
$\varphi: R^n\to R^m$ a homomorphism of free,
finitely generated $R$-modules. Suppose 
that the ideals 
$\langle k\times k\,\text{minors of}\ \varphi
\rangle$ are principal for all $k$. Then
$\J:=\im\varphi$ is free.
\end{coisa}

\begin{proof}
Let 
\\\centerline{$
A=\bsm
a_{11}&\dots&a_{1n}\\
\vdots&\ddots&\vdots\\
a_{m1}&\dots&a_{mn}
\esm$} 
be the $m\times n$  matrix associated to $\varphi$, i.e. the columns of
$A$ generate $\J$.
By hypotheses for $k=1$, the ideal of coefficients of $A$ is
principal:
\\\centerline{$
\ide{a_{11},\dots,a_{ij},\dots,a_{mn}}=
\ide{b} 
$}
We may assume $b\ne0$. 
Let $b_{ij}:=\frac{a_{ij}}{b}$. We may suppose
$a_{11}=b$. Let 
$\J'$ be the module generated by the columns
of 
\\\centerline{
$B=\bsm
1&\dots&b_{1n}\\
\vdots&\ddots&\vdots\\
b_{m1}&\dots&b_{mn}
\esm$.}

\medskip\noindent
Equivalently  (by elementary operations) $\J'$ is generated by the
columns of 
$\bsm
1&0\\
0&B'\\
\esm
$, 
where
\\\centerline{
$B'=\bsm
b_{22}&\dots&b_{2m}\\
\vdots&\ddots&\vdots\\
b_{m2}&\dots&b_{mn}
\esm
$.} 

\medskip\noindent
Applying the inductive hypotheses, we have that
$\im B'$ is free. Thus
$\J'$ is free. Since $R$ is a domain we have
$\J=b\J'\simeq \J'$ and we conclude that $\J$ 
is free.
\end{proof}

\coi{Remark}\label{blowup}
The variety $\mathbb{B}$ is the well known
variety of complete conics (see 
\,\cite{kt},
 \cite{Isqua}). It is constructed
to solve the indeterminacies of the map
$$e_{1}:\p 5\dashrightarrow \check{\p {5}}
$$ sending
a conic to the envelope of its tangent
lines.
Equivalently $e_{1}(A)=\wed2\! A$,
where $A$ is the symmetric
matrix of the conic. We have that $\mathbb{B}$
is equal to the
closure of the graph of   ${e_{1}}$
in $\p 5\times\check{\p
 {5}}$. 
\eco

\coi{\bf A parameter space for foliations with
invariant smooth conic}

Consider  the projective bundle
associated to $\E$ cf.\,\rref{invconics}, and let $q_{2}:\p{}(\E)\to \p N$
denote the projection. 
Then   
\\\centerline{$
\Y=q_{2}(\p{}(\E))\subset\p N
$} 
is a compactification for the parameter space
of foliations   with an invariant smooth conic.
It's not difficult to show that $q_{2}$ is
generically injective, so the  
degree of $\Y$ is given by the top
dimensional Segre class $s_{5}(\E)$. It will be
computed using Bott's formula.

\eco

\begin{coisa}{\bf Theorem.}
Notation as above, let \Y\ be
the compactification for the parameter space
of 1-dimensional foliations of degree d on \p2
 with an invariant smooth conic. Then
the degree of $ \Y$ is given by 
$$\frac{1}{2^55!}(d-1)d(d+1)(d^7+25d^6+231d^5
+795d^4+1856d^3+2468d^2+2256d+768)
$$
and its codimension is equal to $2(d-1)$.
\eco

\begin{proof}

From the definition of $\E$
(cf.\,\ref{invconics}) we see that
$\rk(\E)=d(d+2)$. As $q_{2}$ is generically
injective it is easy to see that $\cod
\Y=2(d-1)$.

To compute the degree  we use Bott's formula,
\be\label{bott}
\int
s_{5}(\E)\cap [\mathbb{B}]=
\sum_{p\in\mathbb{B}^T}
\frac{s_{5}^T(\E_{p})\cap
[p]}{c_{5}^T(T_{p}\mathbb{B})}
\ee
(cf.\,\cite{brion})
where the torus 
$T:=\C^*$ acts on $\mathbb{B}$ with
isolated fixed points. We recall that each
$T$-equivariant Chern class $c_i^T$
 appearing in the
left hand side is simply the $i$-th symmetric 
function on the weights of the fiber of the
 vector bundle at the fixed point.
The action of $T$ on $\B$ is induced by the
action of $T$ on 
$S_{1}$,
given by 
$$t\cdot z_{i}
=t^{w_{i}}z_{i}
$$
for  a suitable choice of weights
$
w_{0},w_{1},w_{2}\in \Z$.
This action induces
 actions on $\p {}(S_{1})$,
on $\p {}(S_{2})=\p 5$ and on $\p
{}(\s{2}\wed 2S_{1})=\check{\p{}}^ 5$ in such a
way that the map
$e_{1}:
\p 5\dashrightarrow
\check{\p 5}
$
(see\,\ref{blowup})
is $T$-equivariant.
Thus we obtain an action of $T$ on
$\mathbb{B}=\text{closure of }\,{\graf{e_{1}}}
\subset\p5\vez\pd5
$.

It is easy to see that if we choose the weights
in such a way  that  all sums
\\\centerline{$w_{i}+w_{j}\ \text{ with } \ 
0\leq i\leq j\leq 2
$} are pairwise distinct,
we obtain the following six isolated fixed points in $\p 5$: 
$$z_{0}^2,z_{0}z_{1},z_{0}z_{2},z_{1}^2,z_{1}z_{2},z_{2}^2.
$$ 

It remains 
to find the fixed point on the fiber of
$\pi:\B\ra\p 5$ over each fixed point in $\p
5$. If $A\in 
\{z_{0}z_{1},z_{0}z_{2},z_{1}z_{2}\}$, 
then $\pi^{-1}(A)$ has just one  point.
For example,
$\pi^{-1}(z_{0}z_{1})=(z_{0}z_{1},\z_{2}^2)$.
Here we put $\z_i$ for the dual basis.

Take $A\in \{z_{0}^2,z_{1}^2,z_{2}^2\}$, say
$A=z_{0}^2$.  Recall that the
exceptional divisor of the blowup is $E:=\p
{}(\N)$ where $\N=\N_{\V{}}\p 5 $ stands for
the normal bundle. We have an explicit
description for $\N$ 
see~\cite[Proposition\,4.4.]{Isqua}.
Notation as in \rref{tau} with $k=1,\,n=2$,
we have
$$\N_{\V{}}\p 5=\OX_{\pd2}(2)\otimes
\s{2}(\T^{\vee}).
$$
The fiber of $\OX_{\pd2}(2)$ (resp. 
$\s{2}(\T^{\vee})$) over $z_{0}$
is the dual space $\ide{\z_{0}^2}$
(resp. $\langle
z_{1}^2,z_{1}z_{2},$ $z_{2}^2\rangle$).
Thus we get
$$
\pi^{-1}(z_{0}^2)=E_{z_{0}^2}=
\p {}(\ide{ \z_{0}^{2}}\otimes\ide{
z_{1}^2,z_{1}z_{2},z_{2}^2}).
$$ 
We see there are three
fixed points in the fiber of  each fixed point $A\in \V{}$. 
For $A=z_{0}^2$ these points are 
\\\centerline{$
\z_{0}^{2}\otimes
z_{1}^2,\z_{0}^{2}\otimes
z_{1}z_{2},\z_{0}^{2}\otimes z_{2}^2.
$}

Summarizing, we have twelve fixed points in $\mathbb{B}$:
\begin{enumerate}
\item $(z_{i}z_{j},\z_{k}^2)$ with  $0\leq{}
i< j\leq2;\,
 k\neq i,j$; these three lie off $E$;
\item $(z_{i}^{2},\z_{i}^{2}\otimes
 z_{j}z_{k})$ with $j,k\neq i;\, j< k$;
\item $(z_{i}^{2},\z_{i}^{2}\otimes
 z_{j}^2)$ with $i\neq j$.
\end{enumerate}

Next we compute the fibers of $\E$ over each
type of fixed point. 
Suppose that   $B\in \mathbb{B}$ is a fixed
point. The strategy is to take
 a curve $B(t)\in\B$ such that 
\[
\lim_{t\rightarrow0}B(t)=B
\] and such that
$A(t):=\pi(B(t))\in\p 5$ is a
curve of smooth conics for $t\ne0$. Therefore
$\E_{B}$ will be obtained as the limit of
$\E_{B(t)}=
\wt\M_{A(t)}$,
\[\lim_{t\rightarrow 0}\wt\M_{A(t)}=\E_{B}.\]     

This enables us to use the  well known
space 
of vector fields of degree $d$
that leave invariant a smooth
conic $\CC=Z(G)$ (see~\cite{Est}). 
This space is  

\medskip
\noindent$
(\spadesuit)\hfill
\{F_{ij}(\frac{\de G}{\de
 z_{i}}\partial_j-\frac{\de G}{\de
 z_{j}}\partial_i)\,\mid\, F_{ij}\in
S_{d-1}\},\hfill$ 

\medskip\noindent
modulo multiples of the radial vector field.

We will adopt the following notation:
for each subset 
$
J:=\{v_{0},\dots,v_{k}\}\subset
\{z_{0},z_{1},z_{2}\}
$
we set 
\\\centerline{$M_{m}(J)
=\{v_{0}^m,v_{0}^{m-1}v_{1},\dots,v_{k}^m\},
$}
the canonical monomial basis of \ $\s{m}(J)$.
We write simply $M_{m}$ for
$M_{m}(\{z_{0},z_{1},z_{2}\})$.

Set 
$\X_{i,j}:=z_{i}\partial_i-z_{j}\partial_j$.
Notice this is a vector of weight 0, since
$t\cdot{}z_i=t^{w_i}z_i$ whereas
$t\cdot{}\de_i=t^{-w_i}\de_i$.

We now describe suitable 1-parameter families
of smooth conics abutting each type of 
fixed point. 

\medskip
\noindent(1)
{$B_{1}=z_{0}z_{1}$}. We take 
$A(t)=z_{0}z_{1}+tz_{2}^2\in\p5$.  
Using the characterization $(\spadesuit)$ we
obtain that the space  $\wt{\M}_{A(t)}$ 
of vector fields leaving $A(t)$ invariant
is given by
$$\{F_{10}(z_{1}\partial_{1}-z_{0}\partial_{0}),
F_{20}(z_{1}\partial_{2}-2tz_{2}\partial_{0}),
F_{21}(z_{0}\partial_{2}-2tz_{2}\partial_{1})\,\mid
F_{ij}\in S_{d-1}\}.$$

Taking limit as $t\rightarrow 0$, we find   
a basis for $\E_{B_{1}}$: 
$$\{F_{1}\X_{0,1},F_{2}\partial_2\,\mid\, F_{1}\in
M_{d-1}\,,\,F_{2}\in M_{d}\setminus\{z_{2}^d\}
\},
$$
Clearly this basis consists of
$T$-eigenvectors.

\medskip
\noindent(2)
{$B_{2}=(z_{0}^{2},\z_{0}^{2}\otimes z_{1}z_{2})$}.
In this case, we take 
$A(t)=z_{0}^2+tz_{1}z_{2}
$.
With the same procedure as above, we obtain the
following 
basis (of eigenvectors) for $\E_{B_{2}}$:
$$
\{F_{1}z_{0}\partial_1,F_{2}z_{0}\partial_2,F_{3}\X_{1,2}\,\mid\, F_{1},F_{2}\in
M_{d-1}\,,\,F_{3}\in M_{d-1}(\{z_{1},z_{2}\})\}
.$$ 

\noindent(3)
{$B_{3}=(z_{0}^{2},\z_{0}^{2}\otimes z_{1}^2)$}.
In this case  a  curve of
smooth conics that approximates $B_{3}$ is  
$A (t)=z_{0}^2+tz_{1}^2+t^2z_{2}^2.
$ 
As before,  we obtain the following basis of
eigenvectors for $\E_{B_{3}}$:
\\\centerline{$
\{F_{1}z_{0}\frac{\de }{\de
 z_{1}},F_{2}\partial_2\,\mid\,
F_{1}\in M_{d-1}\,,\, F_{2}\in
M_{d}\setminus\{z_{2}^d\}\}.
$}

\mn
In order to handle the denominator in Bott's
formula (\ref{bott}) we obtain, for each fixed point $B$, a
base  consisting of eigenvectors of $T_{B}\B$.

If $B\not\in E$ then  $
T_{B}\B \simeq
T_{\pi(B)}\p {}(S_{2})$.
For example, for  $B_{1}=z_{0}z_{1}$ we have
\begin{align}
T_{B_{1}}\B \simeq
\ide{z_{0}z_{1}}^{\vee}\otimes \ide{
z_{0}^2,z_{0}z_{2},z_{1}^2,z_{1}z_{2},z_{2}^2}.
\nonumber
\end{align}

If $B\in E$, then $B=(A,[v])$ with $A\in \V$
and $v\in \N_{A},\,v\ne0$. Now
\begin{equation}
T_{B}\B= T_{A}\V{}\oplus \Hom(\ide{v},
\frac{\N_{A}}{\ide v})  \oplus \ide v\nonumber
\end{equation}
as $\C\us$-modules.
For $B_{2}=(z_{0}^2,\z_{0}^{2}\otimes{}
z_{1}z_{2})$ we have:
\begin{align}
T_{B_{2}}\B =T_{z_{0}^2}\V{}\oplus 
\ide{\z_{0}^{2}\otimes z_{1}z_{2}}^{\vee}
\otimes \ide{\z_{0}^{2}\otimes{}
z_{1}^2,\z_{0}^{2}\otimes z_{2}^2}
\oplus \ide{\z_{0}^{2}\otimes z_{1}z_{2}}
\nonumber
\end{align}
where $T_{z_{0}^2}\V{}= \ide{z_{0}^{2}}^{\vee}
\otimes \ide {z_{0}z_{1},z_{0}z_{2}}$.
Similarly, for 
$B_{3}=(z_{0}^{2},z_{0}^{2\vee}\otimes
z_{1}^2)$ 
we have 
\begin{align}
T_{B_{3}}\B =
T_{z_{0}^2}\V{}\oplus \ide{z_{0}^{2\vee}\otimes
z_{1}^2}^{\vee}\otimes \ide{ z_{0}^{2\vee}\otimes
z_{1}z_{2},z_{0}^{2\vee}\otimes z_{2}^2}
\oplus \ide{ z_{0}^{2\vee}\otimes z_{1}^2}.
\nonumber
\end{align}

The explicit calculation in Bott's formula is
better left for a script (see
\S\,\ref{scripts}). 
Finally we  use Lemma \ref{bound} below 
which enables us to restrict the computation
just for the first sixteen 
values of $d=2,\dots,17$
and then  interpolate the answers obtained.

\end{proof}

\begin{coisa}{\bf Lemma.}\label{bound}
The sum in the right hand side of 
Bott's formula 
$$\int
s_{5}(\E(d))\cap [\mathbb{B}]=\sum_{B\in\mathbb{B}^T}
\frac{s_{5}^T(\E(d)_{B})\cap 
 [B]}{c_{5}^T(T_{B}\mathbb{B})},
$$  
is 
a
combination of $w_{i}'s$ with polynomial
coefficients in $d$ of degree $\leq 15$.  
\end{coisa}

\begin{proof}
For each fixed point $B$ let $\{\xi_{1}(d),\dots,\xi_{m(d)}(d)\}$ denote
the set of weights of $\E(d)_{B}$.
Since  $s_{5}^T(\E(d)_{B})$ is a polynomial in
$\{c_{k}^T(\E(d)_{B})\mid k=1,\dots,5\}$ it's enough to
prove that each $$c_{k}^T(\E(d)_{B})=\sigma_{k}(\xi_{1}(d),\dots,\xi_{m(d)}(d))$$
is a
combination of $w_{i}'s$ with polynomial
coefficients in $d$ of degree $\leq 3k$.

Recalling Newton's identities 
\\\centerline{$
k\sigma_{k}=\sum_{i=1}^{k}(-1)^{i+1}
\sigma_{k-i}p_{i}
$}
where 
\\\centerline{$
p_{k}(\xi_{1}(d),\dots,\xi_{m(d)}(d)):=\sum_{i=1}^{m(d)}
\xi_{i}(d)^k
$}

\mn
we see that it suffices to prove that
$p_{k}(\xi_{1}(d),\dots,\xi_{m(d)}(d))$ 
is a
combination of $w_{i}'s$ with polynomial
coefficients in $d$ of degree $\leq k+2$.

On the other hand, a careful analysis of the
weights appearing in the basis of $\E(d)_{B}$
at each fixed point shows that these weights
can be separated into sets of the form
\\\centerline{
$\{\text{weights of}\,\,M_{e}(J)\}$ or
$\{\text{weights of}\,\,M_{e}(J)\}+w$ 
}

\mn
where $w$
is a (fixed) combination of $w_{i}'s$;
$e=d,d-1$ and $J=\ide{z_{0},z_{1},z_{2}}$ or
$J=\ide{z_{i},z_{j}}$, $i\neq j$.  From this
the reader may  convince her(him)self that
it's enough to prove the following

\begin{description}
\item[Claim] Let $
m=m(d,n):=\binom{d+n}{n}
$
and 
$
\{\xi_{n,1}(d),\dots,\xi_{n,m}(d)\}
$
be the weights associated to the basis
$M_{d}(\{z_{0},\dots,z_{n}\})$ of $\s {d}(\langle z_{0},\dots,z_{n}\rangle)$.
Then 
$$
p^n_{k}(d):=\sum_{i=1}^{m}
\xi_{n,i}(d)^k
$$ 
is  a
combination of $w_{i}'s$ with polynomial
coefficients in $d$ of degree $\leq k+n$.

\end{description}

To prove the claim we proceed by induction on
$n\geq 1$ and on $k\geq 0$. 

For $n=1$, 
\\\centerline{$
M_{d}(\{z_{0},z_{1}\})=\{z_{0}^d,z_{0}^{d-1}z_{1},\dots,z_{0}z_{1}^{d-1},z_{1}^d\}
$}
so that  $m(d,1)=d+1$. We have 
\begin{align}
&p^1_{k}(d)=\sum_{i=1}^{d+1}\xi_{1,i}(d)^k=\sum_{i=0}^{d}(iw_{0}+(d-i)w_{1})^k\nonumber\\
&=\sum_{i=0}^{d}(i(w_{0}-w_{1})+dw_{1})^k=\sum_{i=0}^{d}\sum_{j=1}^k\binom{k}{j}(i(w_{0}-w_{1}))^j(dw_{1})^{k-j}\nonumber\\
&=\sum_{j=1}^k\binom{k}{j}(dw_{1})^{k-j}(w_{0}-w_{1})^j\sum_{i=0}^{d}i^j.\nonumber
\end{align}
The sum $\sum_{i=0}^{d}i^j$ is
polynomial in $d$ of degree $j+1$, therefore $p^1_{k}(d)$ is a polynomial in $d$ of degree $\leq k+1$.

For $k=0$, we have $p_{0}^n(d)=m(d,n)$,
a polynomial in $d$ of
degree $n$.

For the general case, decompose the basis  $M_{d}(\{z_{0},\dots,z_{n}\})$ 
as the union

\smallskip\centerline{$
z_{0}M_{d-1}(\{z_{0},\dots,z_{n}\})\cup z_{1}M_{d-1}(\{z_{1},\dots,z_{n}\})\cup
z_{2}M_{d-1}(\{z_{2},\dots,z_{n}\})\cup\dots
\cup \{z_{n}^d\}.
$}

\mn
Then the weights are:

\medskip\centerline{$
w_{0}+\{\xi_{n,i}(d-1)\}\cup w_{1}+\{\xi_{n-1,i}(d-1)\}\cup
w_{2}+\{\xi_{n-2,i}(d-1)\}\cup\dots \cup
\{dw_{n}\}.
$}

\medskip
Hence we can write
\begin{align}
&p^n_{k}(d)=\sum_{i=1}^{m(d,n)}(\xi_{n,i}(d))^k=
\sum_{i=1}^{m(d-1,n)}(w_{0}+\xi_{n,i}(d-1))^k
~+\nonumber\\
&\sum_{i=1}^{m(d-1,n-1)}(w_{1}+\xi_{n-1,i}(d-1))^k
+\!\sum_{i=1}^{m(d-1,n-2)}(w_{2}+\xi_{n-2,i}(d-1))^k+
\dots+(dw_{n})^k
\nonumber\\
&=\sum_{j=0}^k\binom{k}{j}w_{0}^jp_{k-j}^n(d-1)+\sum_{j=0}^k\binom{k}{j}w_{1}^jp_{k-j}^{n-1}(d-1)
~+\nonumber\\
&
\hskip1cm
\sum_{j=0}^k\binom{k}{j}w_{2}^jp_{k-j}^{n-2}(d-1)
~+\dots+~(dw_{n})^k.\nonumber
\end{align}
By induction  we conclude  that
$p^n_{k}(d)-p^n_{k}(d-1)$ is a polynomial in
$d$
of degree $\leq k+n-1$, and this implies
that $p^n_{k}(d)$ is a polynomial in $d$ 
of degree 
$\leq k+n$.

\end{proof}

\section{Comments}

Using the
 varieties of complete quadrics of any dimension  in
 $\p n$ (see \cite{Isqua}) it is possible to find a compactification of the space of
 dimension one foliations in $\p n$ that leave invariant a smooth
 quadric.  

For example, in the case of  conics in $\p 3$
we obtain the following. 

\begin{coisa}
{\bf Theorem.}
Let $\Y$ denote the closure in $\p N$ of the
variety of dimension one 
foliations in $\p 3$
that have an invariant smooth conic.  The degree of $ \Y$ is given by 
$$\ba c
\frac{4}{8!3^2}(d-1)d\left[207d^{14}+2763d^{13}+15447d^{12}+54395d^{11}+114847d^{10}+\right.
\\\na5\left.
207891d^{9}+
256737d^{8}+225801d^{7}+164937d^{6}+182101d^{5}+38993d^{4}+\right.
\\\na5\left.
316221d^{3}+248856d^{2}
-118908d-332640\right]\ea
$$
and its 
codimension is equal to $4(d-1)$.
\eco

For invariant 
quadric surfaces in $\p 3$ we find:
\begin{coisa}
{\bf Theorem.}
Let $\Y$ denote the closure in $\p N$ of the
variety of dimension one foliations in $\p 3$
that have an invariant quadric.  The degree of
$ \Y$ is given by
$$\ba c
\frac{1}{9!(3!)^9}(d-1)d(d+1)\left[d^{24}+
81d^{23}+3151d^{22}+77949d^{21}+1369333d^{20}
+\right. 
\\\na4\left.18084843d^{19}+185031133d^{18}
+1481854743d^{17}+9251138050d^{16}+
44737976160d^{15}+\right.
\\\na4\left.168507293704d^{14}+503603726976d^{13}
+1212870415960d^{12}+2353394912904d^{11}+\right.
\\\na4\left.3628929239056d^{10}+4249158105672d^9+
3232639214668d^8+413912636928d^7-\right.
\\\na4\left.2874493287072d^6-3885321416832d^5-1115680433472d^4
+4477695012864d^3+\right.
\\\na4\left.8264265366528d^2+8139069775872d+4334215495680\right]\ea$$
and its
codimension is equal to $(d-1)(d+5)$.
\eco


%% file: scripts1.tex
The script below uses \schub\ to compute the
formula for the degree of the variety of
1-dimensional 
foliations of degree $d$ in \p n with some
invariante $k$-plane.
{\small
\begin{verbatim}
with(schubert);
deg:=proc(k,n)local Ec,SrE,F;
grass(n-k,n+1,c);
Ec:=dual((n+1)-Qc);
SdE:=Symm(d,Ec);
F:=SdE*Qc; 
integral(chern(DIM,F));end;
#example
factor(deg(1,2));#1/8*d*(d+3)*(d+2)*(d+1)
\end{verbatim}
}

%% file: scripts2.tex
\subsection{\sc Singular script}
The calculations for Bott's formula uses \sing.
\label{script2}
{\small\baselineskip5pt
\begin{verbatim}
//download from www.mat.ufmg.br/~israel/Projetos/myprocs.ses
//save and load  
<"myprocs.ses"; //more procedures; save and load
//Sum fractions
proc sumfrac(list ab,list cd){
def l= (ideal(ab[1]*cd[2]+ab[2]*cd[1],
ab[2]*cd[2]));def p=gcd(l[1],l[2]);
return(list(l[1]/p,l[2]/p));}

//Subtract fractions 
proc subfrac(list ab,list cd){
def l= (ideal(ab[1]*cd[2]-ab[2]*cd[1],
ab[2]*cd[2]));def p=gcd(l[1],l[2]);
return(list(l[1]/p,l[2]/p));}

//Multiply fractions
proc mulfrac(list ab,list cd){
def l=  ideal(ab[1]*cd[1],
ab[2]*cd[2]);
def p=gcd(l[1],l[2]);
return(list(l[1]/p,l[2]/p));}

//Inverse of a fraction
proc invfrac(list ab){
return(list(ab[2],ab[1]));}

//Divide a fraction
proc divfrac(list ab,list cd){
return(mulfrac(ab,invfrac(cd)));}

//If vs is a poly then return its monomials, 
//else return the degree d monomials 
proc mon(vs,d)
{if(typeof(vs)=="poly"){
def vvs=pol2id(vs);} else {def vvs=vs;}
return(vvs^d);}  
 
//Compute f^k
proc pot(f,k){def p=f;
 for(int i=1;i<=k-1;i++){p=reduce(p*f,mt);}
return(p);}

//Substitute the weights by given values P[i]
proc pesos(H){ //H= poly or list
  list res;
def mm=matrix(vars(1..(nvars(basering))));
def ty@=typeof(H);
if(ty@=="poly"){
  for(int i=1;i<=size(H);i++){
def i_=H[i]; def iv=leadexp(i_); res[i]=(mm*iv)[1,1];}
return(mymapn(ideal(x(0..n)),ideal(P[1..n+1]),res));}
else{if(ty@=="list" and size(H)==2) {
def num,den=H[1..2];
res=pesos(H[1]);def denn=pesos(den);
for(int i=1;i<=size(res);i++) {
res[i]=res[i]-denn[1];}
return(res);}
else{if(ty@=="ideal")
{res=pesos(H[1]);
for(int i=2;i<=size(matrix(H));i++) {
res=res+pesos(H[i]);}
}return(res);}}}

//Symmetric functions in the weights of H
proc cherns(H){
poly sigma= 1;
for(int j=1;j<=min(rk,size(H));j++){
sigma=reduce((1+H[j]*t)*sigma,std(t^(d+1)));}
kill j;def j=coeffs(sigma,t);
while(size(j)<d+1){
j=transpose(concat(transpose(j),[0]));}
return(ideal(j));
}

//Interpolation
proc interpola(lista){
//list(list(a,b),list(c,d),..)
def s=size(lista);
if (s>1) {
for (int j= 2;j<= s;j++){
poly f=1;
for (int i = 1;i<= j;i++){ 
f=f*(t-lista[i][1]); }
if (gcd(f,diff(f,t))==1) {
poly g=0;for(i= 1;i<=j;i++){
def p=f/(t-lista[i][1]);
g =g +lista[i][2] /
subst(p,t,lista[i][1])*p;}
int ii;
for(i=1;i<=s;i++){
if(subst(g,t, lista[i][1] )<>lista[i][2])
{ii++;}}
see(j,ii,s,deg(g));
if (ii==0){ break;}}
else {ERROR("WRONG DATA1");}}
if (deg(g)==s-1){print(g);return("work harder");}
else {
print(myfact(g));
return(g);}}
else {ERROR("WRONG DATA2");}}
//end of procedures

//load dimension of P^n
int n=2;
//dimension of  conics in P^2 
int d=(n+3)*n/2;kill r;
ring r=(0),(x(0..n),t,c(1..d)),dp;
def mt=std(t^(d+1));
//d-Chern class 
def S(d)=summ(seq("c(i)*t^i",1,d));

//d-Segre class
poly p;	for(int i=1;i<=d;i++){
p=p+(-1)^i*pot(S(d),i);}
S(d)=p;def p=coef(S(d),t);
def S(d)=row(p,2);kill p;
S(d)=S(d)[d..1];

//give values to the weights
intvec P=0,1,3;
//forr(3,n+1,"P[i]=random(2,15)");
def xx=ideal(x(0..n));
poly xxs=summ(seq("x(i)",0,n));
list sfs;
sfs[2]=summ(mon(xxs,2));
//Define  lists TP5,TB=tangent of \B,TV,NV=Normal bundle to the Veronese  
//in each fixed point, and find c_{d}(T\B)= product of the weights.
list TP5,TV,NV,TB;
forr(1,size(sfs[2]),"
  def a=sfs[2][i];
  TP5[i]=list(sfs[2]-a,a);
//test if a is of the form Z_{0}^2 or Z_{0}Z_{1}
  def b=radical(a);
//if a is of the form Z_{0}^2
  if(b[1]<>a)
{def p=list(xxs-b[1],b[1]);
  TV[size(TV)+1]=p;
  def q=subfrac(TP5[i],p);
  NV[size(NV)+1]=q;
//For each vector of NV
for(int j=1;j<=size(q[1]);j++){
// Do the sum TV+O_{N}(-1)
def jj=list(q[1][j],q[2]);
def TT=sumfrac(p,jj);
//Obtain TB
TT=sumfrac(TT,divfrac(subfrac(q,jj),jj));
//Save in the list TB the fixed point and the d-Chern class  
TB[size(TB)+1]=list(list(b[1],jj),prod(pesos(TT)));
}}
else
//if a is Z_{0}Z_{1} then TB=TP5
{TB[size(TB)+1]=list(a,prod(pesos(TP5[i])));}
");
list respostas;
//compute the weights of \widetilde{\N} 
//in the fixed points for d=2..17
//to do interpolation.
// dd is the degree of the foliation

for (int dd=2;dd<=17;dd++){
int rk=int((n+1)*binomial(dd+n,n)
	-binomial(dd+1+n,n));
sfs[dd]=summ(mon(xxs,dd));
sfs[dd-1]=summ(mon(xxs,dd-1));
poly NN;
for(int ii=1;ii<=size(TB);ii++) {
def p=TB[ii][1];

//the fixed point, can be list (in E)
// or poly (not in E) if the point is not in E

if(typeof(p)=="poly") {
//convert product in summ, in order to 
//obtain the factors

p=dotprod(leadexp(p),xx);
//obtain the factors 

def a=p[1];def a1=p-a;
def  z=xxs-a-a1;

//Using the basis found  in the text 
//compute the weights of
// the fiber of \widetilde{\N} 

def W=pesos(list(sfs[dd]-z^dd,z))+pesos(sfs[dd-1]);

//symmetric functions in the weights, 
//with given values!
//i.e. c_{i}^T(\widetilde{\N})

W=cherns(W);
//substitute, in s_{5}, c(i) by the 
//values c_{i}^T(\widetilde{\N}) 

def NUM=subst(S(d)[d],c(1),W[1][2]);
for(int i=2;i<=d;i++){NUM=subst(NUM,c(i),W[1][i+1]);}
//Do the summ of the Bott's formula
NN=NN+NUM/TB[ii][2];}
//If the point is in E, p is a list
else{def p=TB[ii]; def a=p[1][1]; def z=p[1][2][1];
//convert prod in summ
z=dotprod(leadexp(z),xx);
//test if a is of the form q_{2} or q_{3}
//if a is of the form q_{3}
if(size(z)==1) {z=leadmonom(z); def w=xxs-a-z; 
//Using the basis found  in the text compute the weights of
// the fiber of \widetilde{\N} 
def W=pesos(list(sfs[dd]-w^dd,w))+pesos(list(a*sfs[dd-1],z));
//symmetric functions in the weights, with given values!

W=cherns(W);

//substitute  c(i) in s_{5} by the values c_{i}^T(\widetilde{\N}) 

def NUM=subst(S(d)[d],c(1),W[1][2]);
for(int
i=2;i<=d;i++){NUM=subst(NUM,c(i),W[1][i+1]);}

//Do the sum in  Bott's formula:
NN=NN+NUM/TB[ii][2];}//matches z=leadmon(z)

//If a is of the type q_{2}:
else{def z1=z[1];z=z[2];
//Using the basis found  in the text compute 
//the weights of the fiber of \widetilde{\N} 

def W = pesos(list(a*sfs[dd-1],z))+
pesos(list(a*sfs[dd-1],z1))+pesos(mon(xxs-a,dd-1));
W=cherns(W);
def NUM=subst(S(d)[d],c(1),W[1][2]);
for(int i=2;i<=d;i++){NUM=subst(NUM,c(i),W[1][i+1]);}
NN=NN+NUM/TB[ii][2];}}}

//Save in respostas the values obtained for each degree

respostas[size(respostas)+1]=list(dd,NN);}

//Interpolate the list  respostas 
interpola(respostas);
\end{verbatim}}

%% file: biblioinv.tex
\enddocument